\documentclass[12pts]{amsart}
\usepackage[english,francais]{babel}
\usepackage[T1]{fontenc}
\usepackage{amsmath,amsfonts}
\usepackage{amsthm}
\usepackage[utf8]{inputenc}
\usepackage{pifont}
\usepackage[all]{xy}
\usepackage{tikz}
\usetikzlibrary{fit,calc,positioning,decorations.pathreplacing,matrix}
\newtheorem{thm}{Theorem}[subsection]
\newtheorem{cor}[thm]{Corollary}
\newtheorem{lem}[thm]{Lemma}
\newtheorem{que}[thm]{Question}
\newtheorem{conj}[thm]{Conjecture}

\newtheorem{rem}[thm]{Remark}

\title{Gaps in the Milnor-Moore spectral sequence and the Hilali-conjecture}

\date{ July 26  2017}


\begin{document}
\titlepage
\author{Youssef Rami}{\let\thefootnote\relax\footnote{{\it Address}: My Ismail University, Departement of Mathematics \& Informatics, B. P. 11 201 Zitoune,  Mekn\`es, Morocco}\let\thefootnote\relax\footnote{{\it Email}: yousfoumadan@gmail.com ; y.rami@fs-umi.ac.ma}}

\maketitle

\selectlanguage{english}








\begin{abstract} 
 Let $X$ be a simply connected  CW-complex, $(\Lambda V, d)$ its minimal Sullivan model, and   $d_k$  the first non-zero homogeneous part of the differential $d$. In this paper, assuming  $(\Lambda V, d_k)$  elliptic, we prove that there isn't any  gap in  the $E_{\infty}$-term of the  Minor-Moore spectral sequence of $X$. Consequently,  we first deduce that in the cohomology $H^*(X,\mathbb{Q})$ there is no $e_0$-gap  and secondly,  we confirm  the  Hilali conjecture when $V = V^{odd}$  or else  when  $k\geq 3$. This last result   generalizes that of \cite{HM08} where $d$ is strictly homogeneous. 
\end{abstract}
\section{Introduction}
The rational  dichotomy theorem  states that a $1$-connected finite type CW-complex $X$ is either $\mathbb{Q}$-elliptic or $\mathbb{Q}$-hyperbolic \cite{F89}. The first ones are characterized by the inequalities $dim(\pi _*(X)\otimes \mathbb{Q})< \infty$ and $dim(H^*(X,\mathbb{Q}))< \infty$. Although they are not generic, they are subject of many work in rational homotopy theory  and  several conjectures are put around them (refer to \cite{FHT01}, Part VI for details).
 Among these we cite the {\it toral-rank conjecture} which is abbreviated by the TRC conjecture and is due to S. Halperin (see \cite{Mun} for its different versions). Recall first that the action of an  $n$-dimensional torus  on $X$ is said {\it almost-free} if all ensuing isotropy groups  are  finite. 
 The largest integer $n\geq 1$, denoted  $rk(X)$,  for which $X$ admits an almost-free $n$-torus action is called the toral-rank of $X$. In \cite{Hal85}, S. Halperin conjectured  the following relation between  this rank and  the  rational cohomology of $X$:
\begin{conj}\label{Tc} (The toral-rank conjecture): If $X$ is a simply-connected finite type CW-complex, then $dim(H^*(X,\mathbb{Q})\geq 2^{rk(X)}$.
\end{conj} 
Refer to \cite{AH}, if $X$ is a $1$-connected  finite CW-complex, then  $rk_0(X)\leq  -\chi_{\pi}$ where $rk_0(X) = rk(X_{\mathbb{Q}})$ and $\chi_{\pi} = \sum_{k\geq 0}(-1)^kdim (\pi_k(X)\otimes\mathbb{Q})$. Therefore, there is a bridge between the TRC  conjecture and the following one posed by   M. R. Hilali in \cite{Hil90} (cf. \cite{HM08} for more details about this bridge).
\begin{conj}\label{Hc} (The H-conjecture): If $X$ is an elliptic simply-connected CW-complex, then $dim(H^*(X,\mathbb{Q})) \geq dim(\pi _*(X)\otimes \mathbb{Q})$.
\end{conj}

The  conjecture (\ref{Hc}) was resolved by different methods for several classes of spaces, as   pure elliptic spaces \cite{Hil90},   hyperelliptic spaces \cite{BFMM12}, elliptic spaces with  formal dimension $\leq 16$ \cite{OT11} and  formal spaces \cite{HM08}.
 These results are obtained after the translation of (\ref{Hc}) in terms of the  minimal Sullivan  model of $X$ as follow. 
   \begin{conj}\label{AvHC} (The algebraic version H-conjecture): If $(\Lambda V,d)$ is an elliptic minimal Sullivan   algebra, then $dimH(\Lambda V,d) \geq dimV$.
 \end{conj}
  Furthermore, refer to  \cite{HM08}, an immediate application of   \cite[Theorem 2.2(B)]{Lup02} gives a positive answer for (\ref{AvHC}) for any minimal Sullivan  model $(\Lambda V,d)$ whose differential is homogeneous of length $l\geq 3$. When $l=2$, by \cite[Corollary 2.6]{Lup02} this is the case if $ker(d_2 : V^{odd}\rightarrow \Lambda V)$ is non zero (that is, if the rational Hurewicz homomorphism is non-zero in some odd degree).
  
  Our goal in this paper is to extend this last result to all minimal Sullivan  algebras $(\Lambda V,d)$ satisfying the condition:
 \begin{eqnarray}\label{cond}
d = \sum_{i\geq k}d_i,\; (k\geq 2)\; \hbox{and}\; (\Lambda V,d_k)\; \hbox{is elliptic}.
 \end{eqnarray} 
   To this end, we first give another formulation of Theorem 2.2(B) in \cite{Lup02}  in terms of gaps in the Milnor-Moore spectral sequence. Indeed, if  $\Lambda V$  is endowed with the  word length filtration $F^p(\Lambda V) = \Lambda ^{\geq p} V$, it is straightforward   to prove that the induced spectral sequence has the form: 
   \begin{eqnarray}\label{algMMss}
   E_k^{p,q} \cong H^{p,q}(\Lambda V,d_k)\Rightarrow H^{p+q}(\Lambda V,d).
   \end{eqnarray} 
   Notice that since $(\Lambda V,d_k)$ is elliptic and hence $dim(V)<\infty$, the filtration $F^p(\Lambda V)$ is bounded and then (\ref{algMMss}) is convergent. Consequently, $dimH(\Lambda V,d)<\infty$ so that $(\Lambda V,d)$ is also elliptic.
When $(\Lambda V,d)$ is the minimal Sullivan model of a simply connected finite type CW-complex $X$, this later is isomorphic from its second term to the Milnor-Moore spectral sequence (cf. \cite{FH82}): 
\begin{eqnarray}\label{MMss}
Ext_{H_*(\Omega X, \mathbb{Q})}^{p,q}(\mathbb{Q},\mathbb{Q})\Rightarrow H^{p+q}(X,\mathbb{Q}).
\end{eqnarray}
    It is then interesting  to recall  that  Lupton's main motivation in \cite{Lup02} was the following question asked by Y. F\'elix   \cite[Question 3.7]{Lup02}:
  \begin{que}\label{Q}
  Can an elliptic CW-complex have $e_0$-gaps in its cohomology?
  \end{que} 
 Here  the cohomology $H^*(X, \mathbb{Q})$ has an $e_0$-gap if it has an element $x$ whose Toomer invariant $e_0(x) = k$,  but does not have any element whose Toomer invariant is $k-1$ (see \S 2 for the definition of $e_0(x)$). So, in other words, \cite[Theorem 2.2 (B)]{Lup02} shows that the cohomology of any elliptic Sullivan  model $(\Lambda V,d)$ with an homogeneous differential  $d$ is without $e_0$-gaps in its cohomology.
 



Now,   $H^+(\Lambda V,d_k) = \oplus _{p\geq 1}H^+_p(\Lambda V, d_k) $, where  cohomology classes in $H^+_p(\Lambda V,d_k)$ are represented by homogeneous cocycles of length $p$ \cite{Lup02}. This implies that $$H^{p,q}(\Lambda V,d_k) = H^{p+q}_p(\Lambda V, d_k).$$   Hence, assuming  that $(\Lambda V,d_k)$ is elliptic, the aforementioned formulation of  \cite[Theorem 2.2 (B)]{Lup02} is as follow:
\\
\centerline{\it the spectral sequence (\ref{algMMss}) has no gaps in its term $E_k^{*,*}$.}

Notice that in  \cite[Proposision 2.7]{KV}, T. Kahl and L. Vandembroucq  showed that  the first term of (\ref{MMss}) has no gaps  and  provided a non rationally elliptic finite CW-complex which presents gaps in  the $E_{\infty}$ term \cite[Corollary 3.3]{KV}.  
Our main result, stated below, proves that this can't happen for some class of rationally  elliptic spaces: 
\begin{thm}
If $X$ is a simply connected finite type  CW-complex whose  Sullivan   model $(\Lambda V,d)$ is such that $(\Lambda V,d_k)$ is  elliptic, then,  at the $E_{\infty}$ term of the Milnor-Moore spectral sequence (\ref{MMss}) of $X$, there can't be any gap. 
\end{thm} 
As a first consequence, we obtain a generalization of Theorem 2.2 (B) in \cite{Lup02} and, by the way, an invalidation to (\ref{Q}) as follows:
\begin{thm}
Any elliptic CW-complexe whose minimal Sullivan model $(\Lambda V,d)$ is such that  $(\Lambda V,d_k)$ is also elliptic has no $e_0$-gap in its cohomology.  
\end{thm}

Returning to the $H$-conjecture (\ref{Hc}), we  state  the second consequence of theorem $1.0.5$ in the following
\begin{thm} 
Let $X$ be a simply connected CW-complex $X$ with a minimal Sullivan   model $(\Lambda V,d)$ such that $(\Lambda V,d_k)$ is elliptic. Then, the H-conjecture holds if 
 $V = V^{odd}$, or else if
 $k\geq 3$.
  \end{thm}
Notice that, under some restrictive hypothesis,  the case where $V=V^{odd}$ was already asserted by A. Amann (\cite{AM}), whom I thank for informing and inspiring me to add such a generalization here. 



   
  {\it  {\bf  Acknowledgements}:   I would like to thank all  members of the MAAT for the beautiful ambiance and atmosphere of work that prevailed within the group. This work
    is based on several discussions exchanged during their monthly seminar.}

\section{Preliminary}   

 Let  $\mathbb K$ be a field of characteristic zero, $V= \oplus _{i=0}^{i=\infty}V^i$  a
graded $\mathbb K$-vector space and  $\Lambda V =
Exterior(V^{odd}) \otimes Symmetric(V^{even})$ ($V^i$ is the subspace of elements of degree $i$). 
A {\it Sullivan algebra} is a free commutative differential graded algebra $(\Lambda V, d)$ with $deg(d) = +1$ and $V$ has  a well-ordered   basis $\{x_{\alpha}\}$  satisfying
$dx_{\alpha }\in \Lambda V_{<\alpha }$ where $V_{<\alpha} = \{v_{\beta} \mid \beta < \alpha \}$. Such algebra is said {\it
minimal} if in addition $deg(x_{\alpha })< deg(x_{\beta })$ implies $\alpha
<\beta $. When  $V^0=\mathbb{Q}$ and  $V^1=0$, this is equivalent to saying that $
d(V)\subseteq \oplus _{i=2}^{i=\infty} \Lambda ^iV  =: \Lambda ^{\geq 2}V$.
 A {\it Sullivan model} for a commutative differential graded algebra  $(A,d)$
 is a quasi-isomorphism $(\Lambda V, d)\stackrel{\simeq } \rightarrow (A,d)$ (morphism inducing an isomorphism in cohomology)  with  source, a Sullivan
algebra. 

When $\mathbb{K} = \mathbb{Q}$ and $X$ is any simply connected CW-complex, the {\it minimal Sullivan  model} of $X$ is by definition any quasi-isomorphism   $(\Lambda V,d)\stackrel{\simeq } \rightarrow A_{PL}(X)$ where $A_{PL}(X)$ stands for  the algebra of polynomial
differential forms associated to $X$ \cite{Sul78}. It is related to rational homotopy groups of $X$ by 
$V^i\cong Hom_{\mathbb Z}(\pi _i(X), \mathbb Q) ;\;\;\; \forall
i\geq 2$  and $X$ is said {\it rationally elliptic} if both $V$ and $H(\Lambda V,d)$ are finite dimensional.     

Now, let 
$(\Lambda V, d)$ be any minimal Sullivan  algebra and  consider the graded algebra morphism 
$p_n:\Lambda V\rightarrow {\Lambda V}/{\Lambda ^{\geq
n+1}V}$  onto the quotient differential graded
algebra obtained by factoring out the differential graded ideal
generated by monomials of length at least $n+1$. {\it The
 Toomer invariant } $e_{ \mathbb K}(\Lambda V,d)$ of
$(\Lambda V, d)$ is the smallest integer $n$ such that $p_n$ induces
an injection in cohomology or $\infty $ if there is no such
integer.
Refer to \cite{FH82}, $e_{\mathbb{K}}(\Lambda V,d)$ is also expressed in terms of the Milnor-Moore spectral sequence (\ref{MMss}) and hence  by isomorphism in terms of (\ref{algMMss})
 by the formula $e_{\mathbb{K}}(\Lambda V, d)= sup\{p \mid E_{\infty}^{p,q}\not = 0\}$
or $\infty$ if such maximum doesn't exists.

Being of a particular interest for our purpose here, we recall that the cohomology $H(\Lambda V,d)$ of  an elliptic  space  with   formal dimension $N = sup\{ p\mid H^p(\Lambda V,d_k)\not =0\}$  satisfies the
 Poincar\'e-duality property in the sens that there are non-degenerate pairings 
 \begin{eqnarray}\label{Pd}
   H^i(\Lambda V,d)\times H^{N-i}(\Lambda V,d)\rightarrow H^{N}(\Lambda V,d)\cong \mathbb{Q}
 \end{eqnarray}
 for $i=1, 2, \ldots N-1$ \cite[Prop. 38.3]{FHT01}.
 
 Consider now $x \in H(\Lambda V,d)$ an arbitrary non zero cohomology class.  G. Lupton defined in \cite{Lup02} its Toomer invariant $e_0(x)$ to be the smallest integer $n$ for which $p_n^*(x)\not = 0$  so that  
  $e_{\mathbb{K}}(\Lambda V,d) = e_0(\omega ) =: e,$  where $\omega $ stands for the fundamental class of $(\Lambda V,d)$.
 
  Moreover, since $H(\Lambda V,d)$ is bi-graded, each one of the pairings above restricts to the following ones
  \begin{eqnarray}\label{Pdbig}
  H^i_p(\Lambda V,d)\times H^{N-i}_{e-p}(\Lambda V,d)\rightarrow H^{N}_e(\Lambda V,d)\cong \mathbb{Q}
  \end{eqnarray}
   for $p = 1, 2, \ldots e-1$ \cite{LM02}.

\section{Proofs of our results:} Denote by $(\Lambda V,d)$ ($d=\sum
_{i\geq k}d_i$, $k\geq 2$) a minimal Sullivan  model of $X$.
We recall in what follow, some general facts (extracted mainly from unpublished notes  by S. Merkulov entitled "Notes on the theory of spectral sequences") about the  spectral sequence of  a filtered complex  essentially to fix  notations and terminology .
Recall first that the general term of such spectral sequence is given by 
$$E_{r}^{p,q} = Z_{r}^{p,q}/Z_{r-1}^{p+1,q-1} +  B_{r-1}^{p,q},$$ where
$$ Z_{r}^{p,q} = \{ x\in [F^p(\Lambda V)]^{p+q} \mid dx\in [F^{p+r}(\Lambda V)]^{p+q+1}\}$$ and 
$$ B_{r-1}^{p,q} = d([F^{p-r+1}(\Lambda V)]^{p+q-1})\cap F^p(\Lambda V) =  d(Z_{r-1}^{p-r+1,q+r-2}).$$
We verify easily that in (\ref{algMMss}), the  $(k-1)^{-th}$-term   coincides with  $(\Lambda V,d_k)$. Therefore its first term $E_k^{p,q} \cong H(\Lambda V,d_k)$. 

Recall also that the differential $\delta _k : E_k^{p,q}\rightarrow E_k^{p+k,q-k+1}$ in $E_k^{*,*}$ is induced from the differential $d$ in $(\Lambda V,d)$ by the formula $\delta _k[v]_k = [dv]_k$, 
 $v$ being any representative in $Z^{p,q}_k$ of the class $[v]_k$ in $E^{p,q}_k$. 
 
 Put $Z(E_k^{p,q}) := Ker(\delta _k)$ and $B(E_k^{p,q}) := Im(\delta _k)$. As $Z_{k+1}^{p,q} + Z_{k-1}^{p+1,q-1}\subseteq Z^{p,q}_k$, a
  natural  monomorphism $$I^{p,q}_k :  {Z_{k+1}^{p,q} + Z_{k-1}^{p+1,q-1}}/{Z_{k-1}^{p+1,q-1} + d(Z_{k-1}^{p-k+1,q+k-2})} \rightarrow E_k^{p,q}$$ 
is given by $I^{p,q}_k(\bar{v}_1 + \bar{v}_2) = \overline{v_1 + v_2} =: [v_1+v_2]_k$. Since $d^2=0$, we have $d(Z_{k+1}^{p,q})\subseteq Z_{k-1}^{p+k+1,q-k}$ and then $d(Z_{k+1}^{p,q} + Z_{k-1}^{p+1,q-1})\subseteq Z_{k-1}^{p+k+1,q-k} + d(Z_{k-1}^{p+1,q-1}) $ which implies that $d(v_1+v_2)\in Z_{k-1}^{p+k+1,q-k} + d(Z_{k-1}^{p+1,q-1})$. Whence the relation  $\delta _k \circ I^{p,q}_k =0$.
Thus, $Im(I^{p,q}_k) \subseteq Z(E_k^{p,q})$. With a little more analysis, one shows that $I^{p,q}_k$ is a surjection and then we have the isomorphism:
\begin{equation}\label{7}
   I^{p,q}_k: {Z_{k+1}^{p,q} + Z_{k-1}^{p+1,q-1}}/{Z_{k-1}^{p+1,q-1} + dZ_{k-1}^{p-k+1,q+k-2}} \stackrel{\cong}{\rightarrow} Z(E_k^{p,q}) 
\end{equation}
Analogous argument gives the proof of the isomorphism: 
\begin{equation}\label{8} 
J_k^{p,q}: {dZ_{k}^{p-k,q+k-1} + Z_{k-1}^{p+1,q-1} }/{Z_{k-1}^{p+1,q-1} + dZ_{k-1}^{p-k+1,q+k-2}} \stackrel{\cong}{\rightarrow}   B(E_k^{p,q}).
\end{equation}
Finally, by induction, from (\ref{7}) and (\ref{8}) one obtains the isomorphisms $E_{r+1}^*\cong H(E_r^*)$ ($r\geq k$).
\subsection{Proof of Theorem 1.0.5} 
\begin{proof}
As mentioned in the introduction, there is no gaps  in the first term of the Milnor-Moore spectral sequence (\ref{algMMss}), that is, referring again to  notations of  \cite{Lup02},  for  any $p=1, \ldots , e$, there exists a non zero class $[\omega _p]\in H_p^{n_p}(\Lambda V,d_k)$. Let $[\omega _0] = 1_{\mathbb{Q}}$ and $H^{0}_{0}(\Lambda V,d_k) = \mathbb{Q}$. So, by (\ref{Pdbig}), to each $[\omega _p]$, it is associated  another non zero class $[\omega _{e-p}]\in H^{N_{e-p}}_{e-p}(\Lambda V,d_k)$ such that $[\omega _e] = [\omega _p]\otimes [\omega _{e-p}]$. 
Recall also that $n_p = min\{ i \mid H^i_p(\Lambda V,d_k) \not = 0\}$ (resp. $N_{e-p} = max\{ i \mid H^i_{e-p}(\Lambda V,d_k) \not = 0\}$) are such that $n_0 =N_0 = 0$, $n_e = N_e = N$ and $n_p + N_{e-p} = N_e$ ($1\leq p \leq e-1$).

In the remainder, we will identify   $H^{p,q}(\Lambda V, d_k)$ and  $E_{k}^{p,q}$, the first term in (\ref{algMMss}),  hence 
 $[\omega _p]\in E^{p,n_p-p}_k$  and $[\omega _{e-p}]\in E^{e-p,N_{e-p}-e+p}_k$. We can then take  $\omega _p\in Z_{k}^{p,n_p-p}$ and $\omega _{e-p}\in Z_{k}^{e-p,N_{e-p}-e + p}$. 
Let us denote $[\omega _p] =: \bar{\omega} _p$ and $[\omega _{e-p}] =: \bar{\omega} _{e-p}$. 
By Theorem 2.2. (C) and Lemma 2.1. in \cite{Lup02} the integers $n_p$
 and $N_{e-p}$
 satisfy the tow equivalent relations:
$$n_2\geq 2n_1,\; n_3\geq n_2 + n_1,\; \ldots , n_{p+1}\geq n_p + n_1,\; \ldots , n_e\geq n_{e-1} + n_1$$
and 
$$N_e = N_{e-1} + n_1,\; N_{e-1}\geq N_{e-2} + n_1,\; \ldots \;  N_{p+1}\geq N_{p} + n_1,\; \ldots , N_1\geq n_1.$$
Since  $H_1(\Lambda V,d_k)\not = 0$ and $V = V^{\geq 2}$, we have $n_1\geq 2$. Therefore, as $k\geq 2$ we have for any  $p= 1, \ldots  ,e-1$, 
 $$(*)\;\;\;\; n_p +1<n_{p+k}<n_{p+(k+1)}< \ldots < n_e$$ and $$(**)\;\;\;\; N_0< \ldots < N_{e-p-(k+1)} < N_{e-p-k}<N_{e-p}-1.$$ Now using the identification above for any $p= 0, 1, \ldots ,e$, we have $$\delta _k : E_k^{p, n_p-p}=H_p^{n_p}(\Lambda V,d_k)\rightarrow E_k^{p+k,n_p-p-k +1}=H_{p+k}^{n_p+1}(\Lambda V,d_k)=0$$ and similarly $$\delta _k : 0 = H_{e-p-k}^{N_{e-p}-1}(\Lambda V,d_k) =  E_{k}^{e-p-k, N_{e-p}-e+p+k-1}\rightarrow E_{k}^{e-p, N_{e-p}-e+p}.$$ 
  It results that
$\forall \; p = 1, \ldots , e-1$, we have necessarily:
\begin{enumerate}
\item[a)]  $\delta _k (\bar{\omega} _p) = \bar{0}$, so that $\bar{\omega} _p \in Z(E_k^{p,n_p-p})$. 
 By the isomorphism (\ref{7}) and due to its  homogeneity of length $p$,  ${\omega} _p$  is obligatory an  element of $Z_{k+1}^{p,n_p-p}$.  Hence $d(\omega _p) \in F^{p+k+1}(\Lambda V)$.
\item[b)]   $\bar{\omega} _{e-p}$ can't be a $\delta _k$-coboundary  i.e. $\bar{\omega} _{e-p}\notin B(E_k^{e-p,N_{e-p}-e+p})$.
\item[c)] $\bar{\omega} _{e}$ is an $\delta _k$-cocycle  that survives to the $\infty $-term $E_{\infty}^{e, N_e-e}$. In particular we have $\omega _e \in Z_{k+1}^{e,N_e-e}$.
\end{enumerate}
Now, since the filtration $F^p(\Lambda V) = \Lambda ^{\geq p}V$  clearly satisfies the relation $$F^p(\Lambda V)\otimes F^q(\Lambda V)\subseteq F^{p+q}(\Lambda V), \; \forall \; p,\; q\geq 0 $$ (i.e. it is one of filtered graded algebras), the  induced spectral sequence  (\ref{algMMss}) is one of graded algebras. It results from  the identification made above that  $\bar{\omega} _e = \overline{{\omega} _p \otimes {\omega} _{e-p}}$ which implies that  $\delta _k(\overline{{\omega} _p \otimes {\omega} _{e-p}}) = \delta _k (\bar{\omega} _e)  = \bar{0}$. Since   $\bar{\omega} _e$ survives to the term $E_{k+1}^{*,*}$, then by homogeneity and the isomorphism (\ref{7}),  $\omega _p \otimes \omega _{e-p}\in Z_{k+1}^{e,N_e-e}$. Whence $d(\omega _p)\otimes \omega _{e-p} \pm \omega _p \otimes d(\omega _{e-p})\in F^{e+k+1}(\Lambda V)$. But $d(\omega _p)\otimes \omega _{e-p}\in F^{e+k+1}$ and then   $\omega _p \otimes d(\omega _{e-p})\in F^{e+k+1}$ also. It results that $d(\omega _{e-p})\in F^{e-p+k+1}$ and then $\omega _{e-p}\in  Z_{k+1}^{e-p,N_{e-p}-e+p}$. Equivalently by identification due to (\ref{7}), we obtain $\delta _k (\bar{\omega} _{e-p}) = \bar{0}$. That is, $\bar{\omega} _{e-p}$ is a $\delta _k$-cocycle, which by $b)$, survives to $E_{k+1}^{*,*}$. Now, using $a)$ and $c)$ below,  we see that  $\bar{\omega} _p$ survives also to  $E_{k+1}^{*,*}$. The relation $[\omega _e] = [\omega _p]\otimes [\omega _{e-p}]$ is therefore  valid in $E_{k+1}^{e, N_e-e}$. We finish the proof by induction using the inequalities $(*)$ and $(**)$.
 \end{proof}
 \subsection{Proof of Theorem 1.0.6}
 \begin{proof}
 By Theorem 1.0.5, we have $E_{\infty}^{p,n_p -p}\not = 0$  ($p = 1, \ldots ,e$). Using  the convergence of (\ref{algMMss}), we have $E_{\infty}^{p,n_p -p} \cong F^p(H^{n_p}(\Lambda V,d))/F^{p+1}(H^{n_p}(\Lambda V,d))$, where $F^p(H^{n_p}(\Lambda V,d)) = Im(H^{n_p}(F^p(\Lambda V),d)\rightarrow H^{n_p}(\Lambda V,d))$.   It results that  $\forall 1\leq p\leq e$, there exists $[\theta _p]\in F^p(H^{n_p}(\Lambda V,d))\backslash F^{p+1}(H^{n_p}(\Lambda V,d))$ so that $e_0([\theta _p]) = p$. Consequently $H(\Lambda V,d)$ has no $e_0$-gap.
 \end{proof}
\subsection{Proof of Theorem 1.0.7.} 
\begin{proof}
 By Theorem 1.0.6, we have  $dimH^{n_p}_p(\Lambda V,d)\geq 1$  for $p = 1, \ldots ,e$. As $dimH^{0}(\Lambda V,d)  = 1$, it follows that $dimH(\Lambda V,d)\geq e$. Now refer to \cite[Theorem 1.]{LM02}, $e =  dimV^{odd}= + (k-2)dimV^{even}$. Hence, if $V = V^{odd}$,  $e = dim(V^{odd}) $, so the inequality holds. Now, if $k \geq 3$, we have immediately $dimH(\Lambda V,d)\geq  dim V$.
\end{proof}
\subsection{Remark}
\begin{enumerate}
\item 
 By \cite[Theorem 2.5]{Lup02}, if  $d=d_k$ is homogeneous and $ker(d_k  : V^{odd}\rightarrow \Lambda V)$ is non zero (i.e. the Hurewicz homomorphism $h: \pi _*(X)\otimes \mathbb{Q}\rightarrow H_*(X;\mathbb{Q})$ is non-zero in some odd degree) then  $dimH^{*}_p(\Lambda V,d_k)\geq 2$ ($1\leq p \leq e-1$). Hence,  $dimH(\Lambda V,d_k)\geq 2e\geq dim(V)$ \cite{HM08}.
  Unfortunately, for $(\Lambda V,d)$ with $d$ non homogeneous and  $k=2$, we can only deduce from Theorem 1.0.5, that  $dimH(\Lambda V,d)\geq e = dimV^{odd}$.  Indeed,  
   we can't be sure that more than one basis element  in $H^{*}_p(\Lambda V,d_2)$ survives to the $E_{\infty}$-term. 

\item The results cited in this paper permit us to ask for a  possible relation between the $H$-conjecture, the  conjecture 3.4 posed by G. Lupton  in \cite{Lup02} and the question (\ref{Q}) asked  by Y. F\'elix.
\end{enumerate}


\end{document}